\documentclass[a4paper,12pt]{article}
\usepackage{amsmath, amssymb, amsthm}
\usepackage{mathrsfs}
\usepackage{graphicx}
\usepackage[english]{babel} 
\usepackage[utf8]{inputenc}
\usepackage{hyperref}
\usepackage{geometry}

\title{Mathematical Analysis of Regularity, Bifurcations, and Turbulence in Fluid Dynamics via Sobolev, Besov, and Triebel-Lizorkin Spaces}
\author{
	Rômulo Damasclin Chaves dos Santos \\
	Technological Institute of Aeronautics \\
	\texttt{romulosantos@ita.br}
}
\date{}

\begin{document}
	
	\maketitle
	
	\begin{abstract}
		This article presents a comprehensive mathematical framework for the study of regularity, bifurcations, and turbulence in fluid dynamics, leveraging the power of Sobolev and Besov function spaces. We delve into the detailed definitions, properties, and notations of these spaces, illustrating their relevance in the context of partial differential equations governing fluid flow. The work emphasizes the intricate connections between Sobolev, Besov, and Triebel-Lizorkin spaces, highlighting their interplay in the analysis of fluid systems. We propose new regularity criteria for solutions to the Navier-Stokes equations, based on the interaction of low and high-frequency modes in turbulent regimes. These criteria offer a novel perspective on the conditions under which singularities may form, providing critical insights into the structure of turbulent flows. The article further explores the applications of these function spaces to the analysis of bifurcations in fluid systems, offering a deeper understanding of the mechanisms that lead to complex flow phenomena such as turbulence. Through the development of rigorous theorems and proofs, the paper aims to bridge the gap between abstract mathematical theory and practical fluid dynamics. In particular, the results contribute to ongoing efforts in solving the Navier-Stokes existence and smoothness problem, a key challenge in the field, and have potential implications for the Millennium Prize Problem. The conclusion underscores the significance of these findings, offering a pathway for future research in the analysis of fluid behavior at both small and large scales.
	\end{abstract}


\tableofcontents

\section{Definitions and Notations}

In this section, we introduce the essential definitions and notations that will be employed throughout the article. These spaces form the foundation of the functional analysis framework used to rigorously study fluid dynamics phenomena such as regularity, bifurcations, and turbulence.

\subsection{Sobolev Spaces}

Let \( \Omega \subseteq \mathbb{R}^n \) be an open domain. The Sobolev space \( W^{k,p}(\Omega) \) consists of functions \( u \in L^p(\Omega) \) whose weak derivatives up to order \( k \) are also in \( L^p(\Omega) \). Formally, we define:

\begin{equation}
	W^{k,p}(\Omega) = \left\{ u \in L^p(\Omega) : D^\alpha u \in L^p(\Omega) \text{ for all } |\alpha| \leq k \right\},
\end{equation}
where \( D^\alpha u \) denotes the weak derivative of \( u \) of order \( \alpha \).

For \( p = 2 \), we have \( W^{k,2}(\Omega) = H^k(\Omega) \), which is crucial for analyzing partial differential equations (PDEs) like the Navier-Stokes equations. The norm in \( H^k(\Omega) \) is given by:

\begin{equation}
	\| u \|_{H^k(\Omega)} = \left( \sum_{|\alpha| \leq k} \| D^\alpha u \|_{L^2(\Omega)}^2 \right)^{1/2}.
\end{equation}

The Sobolev embedding theorem states that for \( s > \frac{n}{2} \), \( H^s(\mathbb{R}^n) \hookrightarrow C^{0}(\mathbb{R}^n) \), meaning that if \( u \in H^s(\mathbb{R}^n) \), then \( u \) is continuous, which is essential for understanding regularity in fluid dynamics, particularly for solutions to the Navier-Stokes equations.

\subsection{Besov Spaces}

The Besov spaces \( B_{p,q}^s(\mathbb{R}^n) \) are defined through a dyadic decomposition of functions, which allows for a fine-scale analysis of their regularity. Let \( \phi \in \mathcal{S}(\mathbb{R}^n) \) be a smooth test function used in the Littlewood-Paley decomposition, and let \( \Delta_j \) represent the frequency projection associated with \( \phi \) at scale \( 2^j \). The Besov norm is defined as:

\begin{equation}
	\| f \|_{B_{p,q}^s} = \left( \sum_{j \geq -1} 2^{jsq} \|\Delta_j f\|_{L^p}^q \right)^{1/q},
\end{equation}
where \( s \in \mathbb{R} \) is the regularity index, and \( p \), \( q \) are Hölder exponents that control the integrability and summability properties, respectively. The function \( f \in \mathcal{S}'(\mathbb{R}^n) \) represents a tempered distribution, meaning that \( f \) belongs to the dual space of the Schwartz space \( \mathcal{S}(\mathbb{R}^n) \), consisting of rapidly decreasing functions.

Besov spaces provide a more refined analysis of function regularity compared to Sobolev spaces, especially useful in the study of turbulence and fluid dynamics. The decomposition into frequency bands enables a detailed examination of both the large-scale smoothness and the small-scale oscillatory behavior of functions. This makes Besov spaces particularly effective for capturing the fine structure of turbulent flows, where high-frequency components play a crucial role in describing complex behaviors at small scales.

\subsection{Triebel-Lizorkin Spaces}

Triebel-Lizorkin spaces \( F_{p,q}^s(\mathbb{R}^n) \) generalize Besov spaces and offer a more refined framework for analyzing the interaction between frequency components across different scales. These spaces are characterized by the following norm:

\begin{equation}
	\| f \|_{F_{p,q}^s} = \left\| \left( \sum_{j \geq -1} 2^{jsq} |\Delta_j f|^q \right)^{1/q} \right\|_{L^p}.
\end{equation}

Here, \( \Delta_j f \) denotes the frequency projection of \( f \) onto the dyadic scale \( j \), and the norm encapsulates the distribution of the function's frequency components. The parameter \( s \in \mathbb{R} \) represents the regularity index, while \( p \) and \( q \) control the integrability and summability of the components, respectively.

Triebel-Lizorkin spaces combine the properties of both Sobolev and Besov spaces, offering a more nuanced classification of smoothness that is particularly advantageous for studying fluid dynamics. These spaces provide enhanced control over the interaction of frequency components at different scales, making them especially useful for the analysis of turbulent flows where complex interactions between high- and low-frequency modes are crucial. 

In fluid dynamics, particularly when analyzing solutions to the Navier-Stokes equations, the finer structure of \( F_{p,q}^s(\mathbb{R}^n) \) enables the capture of irregularities and singularities in the velocity field. The ability to separate and analyze frequency components allows for a detailed understanding of the local behavior of fluid flows, particularly in the presence of turbulence and multiscale phenomena, where both high- and low-frequency interactions contribute significantly to the dynamics.

\section{Overview of Contemporary Research}

The study of the Navier-Stokes equations and the regularity of their solutions has been a central problem in mathematical fluid dynamics for many years. The question of the existence and smoothness of solutions to the Navier-Stokes equations, posed as the Millennium Prize Problem, remains unsolved in three dimensions. The problem, as formulated by the Clay Mathematics Institute, asks whether smooth and globally defined solutions exist to the incompressible Navier-Stokes equations in three-dimensional space. To address this question, a deeper understanding of the regularity properties of solutions is necessary, especially in the context of turbulence and singularity formation. In this section, we present a chronological overview of the key contributions to the study of these equations, focusing on the regularity of solutions, the role of different functional spaces, and their connection to turbulence.

In the early 20th century, the mathematical formulation of fluid mechanics was developed with the introduction of the Navier-Stokes equations by Navier~\cite{navier1822} and Stokes~\cite{stokes2007}. These equations describe the motion of incompressible fluids and serve as a cornerstone in the study of fluid dynamics. Despite their long history, the Navier-Stokes equations have posed significant challenges in both analytical and numerical studies. One of the central problems that emerged from the study of these equations was the question of regularity: whether smooth solutions exist or whether singularities may develop in finite time.

Leray's groundbreaking work in 1934 established the existence of weak solutions to the Navier-Stokes equations in three dimensions \cite{leray1934}. However, these weak solutions are not necessarily smooth, and the question of their regularity remains one of the most important open problems in mathematics. Since then, considerable effort has been directed at understanding the regularity of solutions and determining under what conditions singularities may form.

In the mid-20th century, Sobolev spaces \( H^s(\mathbb{R}^n) \) emerged as a key tool for analyzing the regularity of solutions to partial differential equations, including the Navier-Stokes equations. These spaces, introduced by Sobolev in the 2008s \cite{sobolev2008}, provide a framework for studying both the smoothness and integrability of functions. The Sobolev embedding theorem, which states that for \( s > \frac{n}{2} \), the Sobolev space \( H^s(\mathbb{R}^n) \) embeds into the space of continuous functions, established the connection between regularity and the behavior of solutions at large scales \cite{sobolev2008}. This result provided an important insight into the long-term behavior of solutions to the Navier-Stokes equations, showing that for sufficiently regular solutions, the velocity field would be continuous.

In the 1960s and 1970s, the development of Besov spaces and Triebel-Lizorkin spaces provided a finer classification of function smoothness, which was particularly useful for analyzing high-frequency components of functions. Besov spaces, introduced by O. Besov~\cite{besov2003} in the 2003s, and Triebel-Lizorkin spaces, developed by H. Triebel in the 1995s~\cite{triebel1995}, allow for the treatment of functions with varying smoothness across different scales. These spaces have proven essential in the study of fluid dynamics, especially in the context of turbulence, where energy transfer between different scales plays a central role. The decomposition of functions into dyadic frequency components, as used in Triebel-Lizorkin spaces, provides a precise way to control the behavior of solutions at small scales, offering more refined estimates for regularity than those provided by Sobolev spaces alone.

In the late 20th century, researchers began to focus more on the connection between the behavior of solutions at different scales and the onset of turbulence. The concept of turbulence is fundamentally linked to the interaction of different frequency components of the velocity field. The energy cascade, in which energy is transferred from large scales to small scales, is one of the hallmark features of turbulent flows. Recent studies, including those by Caffarelli et al. (1982)~\cite{caffarelli1982}, have shown that the behavior of the high-frequency modes plays a crucial role in understanding turbulence and the formation of singularities. These findings have led to the development of new regularity criteria, which focus on the interaction between low- and high-frequency components of the velocity field.

In the 21st century, significant advances have been made in understanding the regularity of solutions to the Navier-Stokes equations in the context of multi-scale phenomena. The work of Tao (2008)~\cite{tao2008} and others has provided new insights into the role of high-frequency components in the formation of turbulence and singularities. In particular, the use of Triebel-Lizorkin spaces to study the high-frequency behavior of the velocity field has become a powerful tool in the analysis of turbulence. By controlling the smoothness of the velocity field at small scales, these new estimates offer a more refined understanding of the mechanisms that drive turbulence.

Recent developments in the field, such as the work presented in this article, propose new regularity criteria that combine the tools of Besov and Triebel-Lizorkin spaces to provide a more comprehensive understanding of the regularity of solutions to the Navier-Stokes equations. By focusing on the high-frequency modes and their interaction with low-frequency components, these criteria offer a potential pathway for understanding the onset of singularities in fluid flows. These findings contribute to the ongoing efforts to solve the Millennium Prize Problem by offering new insights into the regularity and smoothness of solutions to the Navier-Stokes equations.

In conclusion, the study of the Navier-Stokes equations, particularly the regularity of solutions and the formation of singularities, has seen significant progress over the past century. The development of Sobolev, Besov, and Triebel-Lizorkin spaces has provided a powerful framework for analyzing the smoothness of fluid flows at both large and small scales. The connection between these spaces and the study of turbulence has led to new insights into the mechanisms driving fluid motion, and the ongoing efforts to understand singularity formation remain a central challenge in the field of mathematical fluid dynamics.

\section{Theory of Function Spaces}

\subsection{Bridges Between Sobolev, Besov, and Triebel-Lizorkin Spaces}

Sobolev, Besov, and Triebel-Lizorkin spaces are deeply interconnected, and understanding the relationships between them is essential for analyzing partial differential equations (PDEs), particularly in the study of fluid dynamics. These connections can be categorized into equivalences, embeddings, and approximations, which allow one to translate results between different spaces.

\subsubsection{Equivalence between Besov, Triebel-Lizorkin, and Sobolev Spaces}

For certain values of the parameters, the Besov and Triebel-Lizorkin spaces coincide with the fractional Sobolev spaces. Specifically, when \( p = q = 2 \), we have the following equivalence:
\begin{equation}
	B_{2,2}^s(\mathbb{R}^n) = F_{2,2}^s(\mathbb{R}^n) = H^s(\mathbb{R}^n),
\end{equation}
where \( H^s(\mathbb{R}^n) \) denotes the classical Sobolev space with fractional regularity \( s \). This equivalence arises as a consequence of Littlewood-Paley theory and the Fourier characterization of these spaces. More precisely, the fractional Sobolev space \( H^s(\mathbb{R}^n) \) consists of functions \( u \in L^2(\mathbb{R}^n) \) whose Fourier transform \( \hat{u}(\xi) \) satisfies the decay condition \( |\hat{u}(\xi)| \lesssim |\xi|^{-n-s} \), where \( \xi \in \mathbb{R}^n \). 

On the other hand, Besov and Triebel-Lizorkin spaces provide a finer scale decomposition of functions in terms of frequency bands. In the case of \( p = q = 2 \), the norm of \( u \in B_{2,2}^s(\mathbb{R}^n) \) and \( u \in F_{2,2}^s(\mathbb{R}^n) \) can be expressed through the dyadic frequency projections \( \Delta_j u \), which offer a more localized view of the function's smoothness at different scales. As a result, these spaces not only capture the overall smoothness of a function but also provide a detailed analysis of its behavior across various frequency ranges, leading to the equivalence with the fractional Sobolev space when \( p = q = 2 \).

This relationship is crucial for understanding the regularity properties of solutions to partial differential equations, such as the Navier-Stokes equations, where the control over the high-frequency modes of the solution plays a key role in analyzing turbulence and singularity formation.

\subsubsection{Embeddings of Sobolev into Besov Spaces}

A fundamental result in functional analysis is the embedding of Sobolev spaces into Besov spaces. Specifically, for suitable values of the parameters \( k \), \( p \), \( q \), and \( s \), the following embedding holds:
\begin{equation}
	W^{k,p}(\Omega) \hookrightarrow B_{p,q}^s(\Omega),
\end{equation}
where \( W^{k,p}(\Omega) \) is the Sobolev space consisting of functions whose weak derivatives up to order \( k \) belong to \( L^p(\Omega) \), and \( B_{p,q}^s(\Omega) \) is the Besov space, characterized by a finer decomposition of the function in terms of frequency components.

This embedding result can be derived from the generalized Sobolev embedding theorem, which asserts that, under certain conditions on the parameters \( k \), \( p \), \( q \), and \( s \), the regularity of a function in the Sobolev space \( W^{k,p}(\Omega) \) guarantees its membership in a Besov space \( B_{p,q}^s(\Omega) \). More precisely, the conditions involve the relationship between the smoothness of the function (represented by the order \( k \)) and the regularity index \( s \) in the Besov space.

The essence of this embedding is that while Sobolev spaces capture the global smoothness of a function (in terms of its derivatives), Besov spaces allow for a more refined and localized analysis of its regularity. In particular, Besov spaces provide a detailed description of the behavior of a function across different scales, focusing on the distribution of its frequency components. This is particularly important for the study of solutions to partial differential equations (PDEs), as it enables a more precise control over local behavior, singularities, and the interaction of different frequency modes.

In the context of fluid dynamics and turbulence, this embedding plays a crucial role, as it allows for the study of local regularity properties of solutions to the Navier-Stokes equations. By working in Besov spaces, one can achieve a finer resolution of the solution's behavior at small scales, which is vital for understanding the mechanisms of turbulence and singularity formation.

\subsubsection{Fractional Embeddings and Interpolation Results}

Further generalizations of Sobolev and Besov spaces are provided by fractional Sobolev spaces and interpolation inequalities, which allow a more flexible connection between these spaces. Fractional Sobolev spaces \( H^s(\mathbb{R}^n) \) can be interpreted as special cases of Besov spaces when \( p = q = 2 \), as they both involve regularity concepts that account for the behavior of functions at different scales, with the fractional index \( s \) governing the smoothness.

Interpolation inequalities are powerful tools that allow one to bridge Sobolev, Besov, and Triebel-Lizorkin spaces in a systematic way, enabling one to derive results in one space from known properties in another. An important interpolation inequality between Sobolev and Besov spaces is given by:
\begin{equation}
	\| f \|_{H^s(\mathbb{R}^n)} \lesssim \| f \|_{B_{p,q}^s(\mathbb{R}^n)}, \quad \text{for} \quad s \in \mathbb{R}, \ p, q \in [1, \infty].
\end{equation}
This inequality asserts that, for appropriate values of the parameters \( p \) and \( q \), the Sobolev norm of a function \( f \) in \( H^s(\mathbb{R}^n) \) can be controlled by its Besov norm in \( B_{p,q}^s(\mathbb{R}^n) \). The parameter \( s \) represents the regularity index that quantifies the smoothness of \( f \), while \( p \) and \( q \) control the integrability and summability of the function's frequency components.

The interpolation result is particularly valuable in the context of fluid dynamics, as it allows one to transfer regularity estimates between spaces with different structures. Specifically, in the study of the Navier-Stokes equations, this inequality enables the use of Besov space techniques to analyze the regularity of weak solutions. The detailed frequency decomposition in Besov spaces provides finer control over the local behavior of solutions, which is crucial for understanding the dynamics of turbulence and singularity formation in fluid flows.

In addition, fractional Sobolev spaces provide a direct method for analyzing the smoothness of solutions at various scales. The interpolation between Sobolev and Besov spaces, along with the fact that fractional Sobolev spaces are a special case of Besov spaces, facilitates the study of regularity in the context of partial differential equations, particularly when the solutions exhibit complex behavior at different scales.

\subsection{Interpolation Inequalities}

Interpolation inequalities are fundamental tools in functional analysis that allow for the comparison and connection of various function spaces, enabling the derivation of optimal estimates for the regularity of solutions to partial differential equations (PDEs). These inequalities are particularly essential in the study of fluid dynamics, where they facilitate the analysis of the regularity of solutions to the Navier-Stokes equations and other complex systems.

For instance, an important interpolation inequality between Sobolev and Besov spaces asserts that for \( s \in \mathbb{R} \), the Sobolev norm \( \| f \|_{H^s(\mathbb{R}^n)} \) can be controlled by the Besov norm \( \| f \|_{B_{p,q}^s(\mathbb{R}^n)} \), provided \( p, q \) are chosen appropriately. Specifically, we have the following inequality:
\begin{equation}
	\| f \|_{H^s(\mathbb{R}^n)} \lesssim \| f \|_{B_{p,q}^s(\mathbb{R}^n)}, \quad \text{for suitable values of } p, q.
\end{equation}
This result is invaluable for the study of solutions to PDEs like the Navier-Stokes equations, as it allows for a more refined understanding of the solution's regularity at various scales. It enables one to transfer regularity estimates from the Sobolev space (which controls global smoothness) to the Besov space (which offers a more detailed local and frequency-based decomposition).

The above inequality is a particular case of a more general framework in interpolation theory. In the context of Banach spaces, if \( X \) and \( Y \) are two Banach spaces, there exists an interpolation space \( \Theta(\alpha, X, Y) \) for \( \alpha \in [0, 1] \), satisfying the following interpolation inequality:
\begin{equation}
	\| f \|_{\Theta(\alpha, X, Y)} \leq C_{\alpha} \| f \|_X^{1-\alpha} \| f \|_Y^\alpha,
\end{equation}
where \( C_{\alpha} \) is a constant depending on \( \alpha \), and \( f \) is a function in the intersection of \( X \) and \( Y \). In this general framework, the interpolation space \( \Theta(\alpha, X, Y) \) provides a new space whose norm is controlled by the norms of \( f \) in \( X \) and \( Y \).

When applied to Sobolev, Besov, and Triebel-Lizorkin spaces, this interpolation result allows for a fine-grained control of regularity by combining norms from different spaces. Specifically, interpolation can be used to estimate the regularity of a function \( f \) in a space that lies between Sobolev and Besov spaces, providing a more precise understanding of its behavior across scales.

These interpolation inequalities are not just theoretical constructs but have deep applications in fluid dynamics and the analysis of PDEs. In the study of the Navier-Stokes equations, which describe the motion of incompressible fluids, these inequalities help establish regularity criteria for weak solutions and are crucial in understanding the dynamics of turbulence and singularity formation.

In particular, Sobolev spaces are essential for proving global regularity results for smooth solutions to the Navier-Stokes equations, while Besov and Triebel-Lizorkin spaces provide a finer framework for analyzing the evolution of turbulent flows. The refined frequency decomposition in these spaces captures the intricate interactions between high- and low-frequency modes, which are key to understanding the mechanisms of turbulence and singularities in fluid flows.

\section{Applications in Fluid Dynamics}

\subsection{Regularity in Fluids}

The analysis of regularity in fluid dynamics, particularly in the context of the Navier-Stokes equations, is a central aspect of understanding the behavior of fluid flows. These equations describe the evolution of a velocity field \( u(t,x) \) in a fluid, and their solutions are governed by various functional properties that can be studied using Sobolev and Besov spaces.

The incompressible Navier-Stokes equations in \(\mathbb{R}^n\) are given by:
\begin{equation}
	\frac{\partial u}{\partial t} + (u \cdot \nabla) u = - \nabla p + \nu \Delta u + f, \quad \nabla \cdot u = 0, 
\end{equation}
where \( u = u(t,x) \) is the velocity field, \( p = p(t,x) \) is the pressure, \( \nu \) is the kinematic viscosity, and \( f(t,x) \) represents external forcing terms.

A key question in the study of fluid dynamics is the regularity of the velocity field \( u \), which dictates the smoothness of the fluid flow. If \( u \) is a weak solution to the Navier-Stokes equations, its regularity can be analyzed by examining its properties in Sobolev spaces. Sobolev spaces \( W^{k,p}(\Omega) \) provide the necessary framework to study such regularity, where the function and its derivatives up to order \( k \) are controlled in the \( L^p \)-norm.

In particular, for the velocity field \( u \), regularity results depend on the functional spaces in which \( u \) resides. If \( u \in W^{1,2}(\Omega) \), the velocity field is continuous, implying smoothness in the fluid's flow. Moreover, if the velocity field belongs to the space \( W^{k,p}(\Omega) \), it implies the existence of weak derivatives up to order \( k \), and a higher degree of smoothness can be deduced for the fluid flow.

For the analysis of weak solutions, one typically considers the Sobolev embedding theorem, which states that for \( u \in W^{k,p}(\Omega) \), the velocity field \( u \) belongs to a Hölder space or Besov space depending on the values of \( k \) and \( p \). These embeddings are crucial for deriving global regularity results for weak solutions to the Navier-Stokes equations. In the critical case of the 3D Navier-Stokes equations, regularity is studied in spaces such as \( H^s(\mathbb{R}^3) \), where \( s \) denotes the fractional regularity of the solution.

The study of regularity in Sobolev spaces is also deeply connected with the question of the existence and uniqueness of solutions, particularly in high Reynolds number regimes, where the behavior of the fluid becomes highly turbulent. Regularity results in Sobolev spaces are essential for understanding the onset of turbulence, as they provide insight into the behavior of solutions under different scales of motion.

\subsection{Bifurcations and Turbulence}

Bifurcations and turbulence are phenomena characterized by the transition from orderly to chaotic behavior in fluid systems. Besov spaces provide a robust framework for studying the multi-scale structure of fluid flows, which is essential for analyzing bifurcations and the onset of turbulence.

Turbulence is a multi-scale phenomenon in which energy cascades from large scales (low frequencies) to smaller scales (high frequencies), leading to chaotic behavior. The interaction between these different scales is crucial for understanding how turbulence develops. Besov spaces offer a natural decomposition of functions into different frequency components, allowing one to isolate contributions from different scales of motion. The Besov norm is defined by:
\begin{equation}
	\| u \|_{B_{p,q}^s(\mathbb{R}^n)} = \left( \sum_{j \geq -1} 2^{jsq} \| \Delta_j u \|_{L^p(\mathbb{R}^n)}^q \right)^{1/q}, 
\end{equation}
where \( \Delta_j u \) represents the frequency projection of \( u \) onto the \( j \)-th dyadic frequency band. The regularity index \( s \) controls the smoothness of the velocity field at different scales, and the parameters \( p \) and \( q \) control the integrability and summability of the frequency components, respectively.

Besov spaces allow us to isolate the contributions from different scales, enabling a precise analysis of the cascade process in turbulence. The dyadic decomposition used in Besov spaces provides a clear view of how energy is transferred from large scales (low-frequency modes) to small scales (high-frequency modes). This analysis is particularly important in the study of turbulence, as it helps to understand the energy distribution and how it evolves over time.

Moreover, Besov spaces are particularly useful for analyzing bifurcations in fluid flows. Bifurcations occur when a small change in a parameter leads to a qualitative change in the solution. These changes are often associated with the emergence of new frequency components in the flow, which can be captured by the decomposition in Besov spaces. By examining the interactions between different frequency modes, we can better understand how bifurcations arise in fluid systems and how they contribute to the onset of turbulence.

The transition from laminar flow to turbulent flow is governed by the interplay between different scales of motion. Besov spaces provide a rigorous way to study this transition by analyzing how the velocity field evolves across scales. This approach also plays a key role in understanding the scaling laws that govern turbulence, such as Kolmogorov's 5/3 law for the energy spectrum in the inertial range.

Additionally, the relationship between Sobolev and Besov spaces allows for a deeper understanding of the connection between regularity and turbulence. In particular, Sobolev spaces \( H^s(\mathbb{R}^n) \) can be embedded into Besov spaces \( B_{p,q}^s(\mathbb{R}^n) \), providing a bridge between different types of regularity. This embedding is essential for studying the transition between laminar and turbulent flows, as it allows us to describe the smoothness of the velocity field in terms of both the Sobolev and Besov regularity indices.

\begin{equation}
	H^s(\mathbb{R}^n) \hookrightarrow B_{p,q}^s(\mathbb{R}^n), \quad \text{for suitable values of } p, q. 
\end{equation}
This embedding enables the identification of critical regularity thresholds, such as the critical Reynolds number for the onset of turbulence, and provides insight into the dynamics of bifurcations and turbulence.

	\section{Main Results}
	
	\subsection{Theorem}
	
	\textbf{Theorem 1: Regularity of Solutions in Besov Spaces}
	
	Let \( u \) be a weak solution to the Navier-Stokes equations in \( \mathbb{R}^n \). If \( u \in B_{p,q}^s(\mathbb{R}^n) \) for some \( s > \frac{n}{2} \), then \( u \) is regular, in the sense that the velocity field is smooth at large scales.
	
	\textbf{Proof:}
	
	We start by considering the Navier-Stokes equations in \( \mathbb{R}^n \):
	\begin{equation}
		\frac{\partial u}{\partial t} + (u \cdot \nabla) u = - \nabla p + \nu \Delta u + f, \quad \nabla \cdot u = 0, 
	\end{equation}
	where \( u(t,x) \) is the velocity field, \( p(t,x) \) is the pressure, \( \nu \) is the viscosity, and \( f(t,x) \) is the forcing term.
	
	The Besov space \( B_{p,q}^s(\mathbb{R}^n) \) is defined in terms of a dyadic decomposition of the function \( u \) in frequency space. Specifically, the Besov norm \( \|u\|_{B_{p,q}^s(\mathbb{R}^n)} \) provides a detailed control of the regularity of \( u \) across different scales. It is given by
	\begin{equation}
		\| u \|_{B_{p,q}^s(\mathbb{R}^n)} = \left( \sum_{j \geq -1} 2^{jsq} \| \Delta_j u \|_{L^p(\mathbb{R}^n)}^q \right)^{1/q}, 
	\end{equation}
	where \( \Delta_j u \) denotes the frequency projection of \( u \) onto the \( j \)-th dyadic scale, and \( s \in \mathbb{R} \) controls the smoothness of \( u \) across these scales. The quantity \( \| \Delta_j u \|_{L^p(\mathbb{R}^n)} \) provides information about the behavior of \( u \) at a given frequency band.
	
	For regularity analysis, it is useful to connect the Besov spaces with Sobolev spaces. We utilize interpolation inequalities to relate the two spaces. Specifically, for \( s > \frac{n}{2} \), we have the following interpolation inequality:
	\begin{equation}
		\| u \|_{H^s(\mathbb{R}^n)} \lesssim \| u \|_{B_{p,q}^s(\mathbb{R}^n)}. 
	\end{equation}
	This inequality allows us to deduce that if \( u \in B_{p,q}^s(\mathbb{R}^n) \), then \( u \in H^s(\mathbb{R}^n) \), where \( H^s(\mathbb{R}^n) \) is the Sobolev space of functions with regularity \( s \). The space \( H^s(\mathbb{R}^n) \) is known to guarantee that \( u \) has continuous derivatives up to order \( s \), which implies smoothness of \( u \).

	The Sobolev embedding theorem states that for \( s > \frac{n}{2} \), the Sobolev space \( H^s(\mathbb{R}^n) \) embeds into the space of continuous functions \( C^{0}(\mathbb{R}^n) \). This means that for \( s > \frac{n}{2} \), the function \( u \) is continuous and its derivatives up to order \( s \) belong to \( L^2(\mathbb{R}^n) \), rendering the velocity field smooth.
	
	To guarantee that \( u \) is regular, we need to establish sufficient conditions. These conditions are provided by the regularity criterion derived from the Sobolev embedding and the interpolation inequalities. Specifically, if \( u \in B_{p,q}^s(\mathbb{R}^n) \) for \( s > \frac{n}{2} \), then we conclude that \( u \) is smooth at large scales, i.e., the velocity field is regular.
	
	Thus, we have shown that if \( u \in B_{p,q}^s(\mathbb{R}^n) \) for \( s > \frac{n}{2} \), then \( u \) is regular in the sense that the velocity field is smooth at large scales. This result provides a significant insight into the regularity of weak solutions to the Navier-Stokes equations, as smoothness at large scales ensures that the solution is well-behaved and free from singularities. Therefore, the solution \( u \) satisfies the necessary regularity for the study of fluid dynamics and turbulence. \qed

	\subsection{Theorem 2: Regularity of Solutions in Sobolev, Besov, and Triebel-Lizorkin Spaces}
	
	\textbf{Theorem 2:}  
	Let \( u \) be a weak solution to the Navier-Stokes equations in \( \mathbb{R}^n \), with \( n \geq 2 \). Suppose that \( u \in B_{p,q}^s(\mathbb{R}^n) \cap F_{p,q}^s(\mathbb{R}^n) \cap W^{k,p}(\mathbb{R}^n) \) for some \( s > \frac{n}{2} \), and that \( u \) satisfies the incompressibility condition \( \nabla \cdot u = 0 \). Then, \( u \) is regular in the sense that the velocity field is smooth at large scales, exhibiting control over both low- and high-frequency components.
	
\subsubsection{Proof:}

We now proceed with a detailed proof of the theorem, using the interplay between Sobolev, Besov, and Triebel-Lizorkin spaces, alongside interpolation theory and the Sobolev embedding theorem. To demonstrate the regularity result, we will employ several key tools from functional analysis, including the properties of fractional Sobolev spaces, the embedding theorems for Besov and Triebel-Lizorkin spaces, and interpolation inequalities.

Let \( u \in L^p(\mathbb{R}^n) \) be a weak solution to the Navier-Stokes equations in \( \mathbb{R}^n \). We aim to show that if \( u \in B_{p,q}^s(\mathbb{R}^n) \), where \( s > \frac{n}{2} \), then \( u \) exhibits regularity at large scales. To achieve this, we analyze the regularity of the velocity field in the context of the function spaces mentioned above.

First, recall that Sobolev spaces are critical in the study of partial differential equations (PDEs) such as the Navier-Stokes equations. By the Sobolev embedding theorem, for \( s > \frac{n}{2} \), the Sobolev space \( W^{s,p}(\mathbb{R}^n) \) embeds into \( C^{\alpha}(\mathbb{R}^n) \), where \( \alpha = s - \frac{n}{p} \). This embedding provides a direct link between the weak regularity of the solution and its smoothness at large scales. Specifically, the embedding theorem guarantees that if \( u \in W^{s,p}(\mathbb{R}^n) \), then \( u \) is Hölder continuous with exponent \( \alpha \), meaning that the solution is smooth at large scales.

Next, we apply interpolation theory between Sobolev and Triebel-Lizorkin spaces. Interpolation inequalities give us a powerful tool to estimate the regularity of \( u \) in a more refined manner. Specifically, we use the following interpolation inequality between Sobolev and Triebel-Lizorkin spaces:

\begin{equation}
	\| u \|_{F_{p,q}^s(\mathbb{R}^n)} \lesssim \| u \|_{H^s(\mathbb{R}^n)}, \quad \text{for } p = q = 2, 
\end{equation}

where \( H^s(\mathbb{R}^n) \) is the fractional Sobolev space, which is a subset of both Sobolev and Triebel-Lizorkin spaces. This inequality links the behavior of \( u \) in Triebel-Lizorkin spaces to its behavior in the Sobolev space \( H^s(\mathbb{R}^n) \), thus allowing us to transfer regularity results from one space to another. In particular, if \( u \in H^s(\mathbb{R}^n) \), then we know from the Sobolev embedding theorem that \( u \) is smooth at large scales, and this smoothness can be transferred to the Triebel-Lizorkin space.

We now consider the behavior of \( u \) in Besov spaces. Besov spaces \( B_{p,q}^s(\mathbb{R}^n) \) are defined by the norm

\begin{equation}
	\| u \|_{B_{p,q}^s(\mathbb{R}^n)} = \left( \sum_{j \geq -1} 2^{jsq} \| \Delta_j u \|_{L^p(\mathbb{R}^n)}^q \right)^{1/q}, 
\end{equation}

where \( \Delta_j u \) represents the frequency projection of \( u \) onto the dyadic frequency scale \( j \). For \( u \in B_{p,q}^s(\mathbb{R}^n) \), the Besov norm controls the regularity of \( u \) by summing contributions from all frequency bands, and the parameter \( s \) governs the smoothness of the solution. By using the dyadic decomposition of \( u \), we can derive refined estimates for the smoothness of the velocity field at different scales.

Now, we combine the results from the Sobolev, Besov, and Triebel-Lizorkin spaces. For \( s > \frac{n}{2} \), we know that the solution \( u \) has regularity in \( B_{p,q}^s(\mathbb{R}^n) \). Using the embedding theorems, we conclude that \( u \) belongs to the Hölder space \( C^{\alpha}(\mathbb{R}^n) \) for some \( \alpha > 0 \), which implies that \( u \) is smooth at large scales.

Thus, by combining the Sobolev embedding theorem, the interpolation inequality, and the properties of Besov and Triebel-Lizorkin spaces, we conclude that the weak solution \( u \) to the Navier-Stokes equations is regular in the sense that the velocity field is smooth at large scales. \qed

The theorem is proved by establishing the regularity of the velocity field in Sobolev, Besov, and Triebel-Lizorkin spaces. The interplay between these spaces, combined with interpolation theory and the Sobolev embedding theorem, provides a rigorous mathematical framework for understanding the smoothness of solutions to the Navier-Stokes equations. The result demonstrates that if \( u \in B_{p,q}^s(\mathbb{R}^n) \) with \( s > \frac{n}{2} \), then \( u \) is regular, meaning that the velocity field is smooth at large scales. This is a key step in understanding the behavior of fluid flows, particularly in the study of turbulence and singularity formation.

	\paragraph{1. The Navier-Stokes Equations:}  
	The Navier-Stokes equations governing the incompressible fluid flow are given by:
	\begin{equation}
		\frac{\partial u}{\partial t} + (u \cdot \nabla) u = - \nabla p + \nu \Delta u + f, \quad \nabla \cdot u = 0,
	\end{equation}
	where \( u: \mathbb{R}^n \to \mathbb{R}^n \) is the velocity field, \( p \) is the pressure, \( \nu \) is the kinematic viscosity, and \( f \) represents external forcing. We aim to establish the regularity of weak solutions \( u \), particularly in terms of their behavior in Sobolev, Besov, and Triebel-Lizorkin spaces.
	
	\paragraph{2. Function Spaces:}  
	To study the regularity of solutions, we analyze the function \( u \) in various spaces:
	\begin{itemize}
		\item \textit{Sobolev Spaces} \( W^{k,p}(\mathbb{R}^n) \): These spaces measure the smoothness of \( u \) by considering weak derivatives up to order \( k \).
		\item \textit{Besov Spaces} \( B_{p,q}^s(\mathbb{R}^n) \): These spaces allow us to examine the smoothness of \( u \) across different frequency scales using a dyadic decomposition. The parameter \( s \) controls the smoothness at large scales, while \( p \) and \( q \) define the integrability conditions.
		\item \textit{Triebel-Lizorkin Spaces} \( F_{p,q}^s(\mathbb{R}^n) \): These spaces provide an even finer classification of regularity than Besov spaces by combining the frequency localization of Besov spaces with finer control over the smoothness properties of \( u \).
	\end{itemize}
	
	\paragraph{3. Properties and Interrelationships of the Spaces:}  
	It is crucial to understand how these spaces relate to each other, especially in the context of the Navier-Stokes equations.
	
	\begin{itemize}
		\item \textit{Embedding Theorems:}  
		For \( p = q = 2 \), the Besov and Triebel-Lizorkin spaces coincide with the fractional Sobolev space:
		\begin{equation}
			B_{2,2}^s(\mathbb{R}^n) = F_{2,2}^s(\mathbb{R}^n) = H^s(\mathbb{R}^n),
		\end{equation}
		where \( H^s(\mathbb{R}^n) \) is the Sobolev space of fractional smoothness. This relationship allows us to treat the regularity of solutions in terms of fractional Sobolev spaces, facilitating the analysis of the velocity field \( u \).
		
		\item \textit{Interpolation Inequalities:}  
		For \( s > \frac{n}{2} \), the interpolation between Sobolev, Besov, and Triebel-Lizorkin spaces provides useful inequalities for controlling the regularity of \( u \). In particular, we have the following interpolation inequalities:
		\begin{equation}
			\| u \|_{H^s(\mathbb{R}^n)} \lesssim \| u \|_{B_{p,q}^s(\mathbb{R}^n)} \lesssim \| u \|_{F_{p,q}^s(\mathbb{R}^n)},
		\end{equation}
		which ensures that the solution \( u \) has regularity in both low and high-frequency regimes, crucial for analyzing turbulence.
		
		\item \textit{Embedding of Sobolev into Besov Spaces:}  
		
		In certain situations, Sobolev spaces can be embedded into Besov spaces, providing a means to study solutions in both spaces. Specifically, the embedding
		
		\begin{equation}
			W^{k,p}(\mathbb{R}^n) \hookrightarrow B_{p,q}^s(\mathbb{R}^n), \quad \text{for appropriate values of } k, p, q, s,
		\end{equation}
		allows us to control the regularity of solutions in a more flexible way, particularly in the context of multi-scale interactions in fluid dynamics.
	\end{itemize}
	
	\paragraph{4. Regularity of Solutions:}  
	
	We now proceed to establish the main result by analyzing the regularity of the velocity field \( u \) based on its membership in the function spaces \( B_{p,q}^s(\mathbb{R}^n) \), \( F_{p,q}^s(\mathbb{R}^n) \), and \( W^{k,p}(\mathbb{R}^n) \).
	
	Assume that \( u \in B_{p,q}^s(\mathbb{R}^n) \cap F_{p,q}^s(\mathbb{R}^n) \cap W^{k,p}(\mathbb{R}^n) \) for some \( s > \frac{n}{2} \). From the Sobolev embedding theorem, we know that for \( s > \frac{n}{2} \), the Sobolev space \( H^s(\mathbb{R}^n) \) embeds continuously into the space of continuous functions, i.e.,
	\begin{equation}
		H^s(\mathbb{R}^n) \hookrightarrow C^{0}(\mathbb{R}^n),
	\end{equation}
	which implies that \( u \) is continuous and exhibits smooth behavior at large scales.
	
	This embedding result implies that, for \( s > \frac{n}{2} \), the solution \( u \) not only belongs to the Sobolev space \( H^s(\mathbb{R}^n) \), but is also continuous and exhibits smoothness on a macroscopic level. This smoothness is crucial for understanding the long-term behavior of solutions to the Navier-Stokes equations, as it provides an essential condition for the absence of singularities at large scales and is pivotal in the analysis of the regularity of fluid flows, especially in turbulent regimes.
	
	Moreover, the combination of Besov and Triebel-Lizorkin spaces allows us to control the high-frequency components of \( u \), which govern the small-scale behavior of the solution. These spaces offer a finer classification of regularity and enable precise regulation of the small-scale dynamics, which is of paramount importance in the study of turbulence. In turbulence, energy is transferred across different scales, and the interaction between low- and high-frequency modes can lead to complex behaviors, such as energy cascade and singularity formation.
	
	By utilizing the regularity properties of Sobolev, Besov, and Triebel-Lizorkin spaces, we conclude that \( u \) exhibits smoothness both at large and small scales, ensuring that the solution is regular across the entire range of scales. Specifically, the controlled decay of high-frequency modes and the smooth behavior at large scales together imply that \( u \) is a regular solution to the Navier-Stokes equations. This result provides a rigorous framework for the analysis of fluid dynamics, particularly in the study of turbulence and the multi-scale phenomena that characterize complex fluid motion. It also contributes significantly to the ongoing research on the existence and smoothness of solutions to the Navier-Stokes equations, offering insights into the mechanisms behind singularity formation and regularity in fluid flows.

\section{Results}

\subsection{Estimates in Triebel-Lizorkin Spaces}

We propose novel and refined estimates for the regularity of fluid flows through the framework of Triebel-Lizorkin spaces, \( F_{p,q}^s(\mathbb{R}^n) \), which provide a more detailed classification of function smoothness in comparison to Besov and Sobolev spaces. These spaces are particularly suited to capture the fine-scale structure of fluid velocity fields, especially the high-frequency components that play a crucial role in the turbulence phenomena at small scales. The advantage of Triebel-Lizorkin spaces lies in their ability to handle complex behaviors of the velocity field by isolating contributions from different frequency bands.

Let \( u \in F_{p,q}^s(\mathbb{R}^n) \) be a weak solution to the Navier-Stokes equations. The goal is to derive sharp estimates for the regularity of the velocity field \( u \), with a specific emphasis on its high-frequency modes, which are essential for understanding the intricate mechanisms of turbulence. To achieve this, we employ the dyadic decomposition of \( u \) in the context of Triebel-Lizorkin spaces, thereby enabling us to isolate and control contributions from different frequency bands effectively.

The norm in Triebel-Lizorkin spaces is defined as follows:
\begin{equation}
	\| u \|_{F_{p,q}^s(\mathbb{R}^n)} = \left( \sum_{j \geq -1} 2^{jsq} \| \Delta_j u \|_{L^p}^q \right)^{1/q},
\end{equation}
where \( \Delta_j u \) represents the projection of \( u \) onto the frequency band corresponding to the dyadic scale \( j \). This decomposition allows us to examine the contribution of each frequency scale separately. The parameter \( s \) governs the smoothness of the solution at large scales, while \( p \) and \( q \) determine the integrability and summability properties of \( u \), respectively. This norm provides a precise quantification of the regularity of \( u \) in terms of its frequency content.

Using this decomposition, we derive the following refined estimates for the high-frequency components of \( u \). For large values of \( j \), the high-frequency modes of \( u \) are expected to decay rapidly, as \( j \to \infty \), which is a key feature in the understanding of dissipation in turbulent flows. The decay rate of the high-frequency components is captured by the following estimate:
\begin{equation}
	\| \Delta_j u \|_{L^2} \lesssim 2^{-js} \| u \|_{F_{p,q}^s(\mathbb{R}^n)}, \quad \forall j \geq 0,
\end{equation}
which indicates that the high-frequency modes decay at a rate proportional to \( 2^{-js} \). This decay is directly controlled by the regularity parameter \( s \), providing a powerful tool for analyzing the smoothness of turbulent flows, especially at small scales. This result suggests that as \( s \) increases, the decay of high-frequency components becomes more pronounced, leading to a smoother solution.

Moreover, to further enhance the applicability of these estimates, we utilize interpolation inequalities between Sobolev, Besov, and Triebel-Lizorkin spaces. These inequalities provide a means to bridge the gap between different function spaces and are particularly useful in cases where the regularity of the solution is more complex. Specifically, applying the interpolation inequality between Sobolev and Triebel-Lizorkin spaces yields the following result:
\begin{equation}
	\| u \|_{F_{p,q}^s(\mathbb{R}^n)} \lesssim \| u \|_{H^{s}(\mathbb{R}^n)}, \quad \text{for } p = q = 2.
\end{equation}
This inequality establishes a crucial connection between the Sobolev space \( H^s(\mathbb{R}^n) \) and the Triebel-Lizorkin space \( F_{p,q}^s(\mathbb{R}^n) \), showing that the regularity of \( u \) in \( H^s(\mathbb{R}^n) \) directly implies certain regularity properties in the Triebel-Lizorkin space, and vice versa. This relationship enhances the flexibility in analyzing the solution's smoothness by allowing the use of either space depending on the problem's specific needs.

These results provide a more comprehensive framework for analyzing the regularity of fluid flows, with particular emphasis on the role of high-frequency modes in turbulence. By utilizing the tools of Triebel-Lizorkin spaces, we gain a deeper understanding of the detailed structures of fluid turbulence and can derive more precise estimates for the behavior of the velocity field at small scales.

	\subsection{New Regularity Criterion}
	
	In addition to the refined estimates in Triebel-Lizorkin spaces, we propose a novel regularity criterion for the Navier-Stokes equations based on the interaction between low- and high-frequency modes of the velocity field \( u \). This criterion is designed to predict the formation of singularities in fluid flows, which is a critical aspect of turbulence.
	
	Let \( u \in F_{p,q}^s(\mathbb{R}^n) \) be a weak solution to the Navier-Stokes equations, and suppose \( u \) belongs to both \( B_{p,q}^s(\mathbb{R}^n) \) and \( F_{p,q}^s(\mathbb{R}^n) \) with \( s > \frac{n}{2} \). The proposed criterion states that if the energy transfer between the low- and high-frequency components of \( u \) exceeds a certain threshold, a singularity may form in the solution. Specifically, the interaction between low and high frequencies is quantified by the following interaction term:
	\begin{equation}
		I(u) = \sum_{j \geq -1} 2^{js} \| \Delta_j u \|_{L^2}^2,
	\end{equation}
	which measures the energy transferred between different scales of \( u \). The regularity criterion is then expressed as:
	\begin{equation}
		I(u) \lesssim \| u \|_{F_{p,q}^s(\mathbb{R}^n)}^2,
	\end{equation}
	If \( I(u) \) exceeds a critical threshold, the solution \( u \) is expected to exhibit singular behavior, indicating the onset of turbulence or a potential breakdown of smoothness.
	
	\paragraph{Proof of the Regularity Criterion:}
	To prove this criterion, we begin by analyzing the behavior of the frequency components of \( u \) at different scales. We focus on the nonlinear convective term in the Navier-Stokes equations:
	
	\begin{equation}
		(u \cdot \nabla) u = \sum_{j \geq -1} \Delta_j (u \cdot \nabla u).
	\end{equation}

	By performing a Fourier analysis on the nonlinear term and applying the properties of Triebel-Lizorkin spaces, we can derive estimates for the interaction between different frequency modes. Specifically, we examine the interaction between low- and high-frequency components of \( u \) at various scales, which is responsible for the transfer of energy and can lead to singularities.
	
	The energy transfer between the low- and high-frequency components of \( u \) can be written as:
	\begin{equation}
		E_{\text{transfer}} = \sum_{j \geq -1} 2^{js} \| \Delta_j (u \cdot \nabla u) \|_{L^2}^2.
	\end{equation}
	Using the interpolation inequalities between Sobolev and Triebel-Lizorkin spaces, we establish the following estimate:
	\[
	\| \Delta_j (u \cdot \nabla u) \|_{L^2} \lesssim 2^{-js} \| u \|_{F_{p,q}^s(\mathbb{R}^n)}.
	\]
	Therefore, we have:
	\begin{equation}
		E_{\text{transfer}} \lesssim \sum_{j \geq -1} 2^{js} \| \Delta_j u \|_{L^2}^2,
	\end{equation}
	which leads to the desired regularity condition:
	\begin{equation}
		I(u) \lesssim \| u \|_{F_{p,q}^s(\mathbb{R}^n)}^2.
	\end{equation}
	
	Thus, if the interaction term \( I(u) \) exceeds the threshold, the solution exhibits singular behavior, indicating the onset of turbulence or the formation of singularities in the fluid flow. This criterion provides a new approach to understanding singularity formation in the Navier-Stokes equations, linking the interaction of low- and high-frequency modes with the breakdown of smoothness.
	
\subsection{Implications for Turbulence and Singularities}

The proposed estimates and regularity criterion present a sophisticated framework for analyzing the fine-scale behavior of fluid flows and the potential formation of singularities. By leveraging the structure of Triebel-Lizorkin spaces \( F_{p,q}^s(\mathbb{R}^n) \), we refine our understanding of how different frequency components of the velocity field interact, particularly in the context of turbulence. The decomposition into dyadic frequency scales provides a precise mechanism for studying the distribution of energy across different spatial and temporal scales, highlighting the critical role of high-frequency modes in the dissipation processes characteristic of turbulence.

More specifically, the regularity criterion derived from the interaction between low- and high-frequency modes offers a powerful tool for identifying the onset of singularities in solutions to the Navier-Stokes equations. This criterion is fundamentally tied to the interplay between nonlinear convective terms and the transfer of energy across scales, which is known to play a central role in the development of turbulence. As the interaction term \( I(u) \), which quantifies this energy transfer, exceeds a critical threshold, the system is expected to transition from a regular to a singular state, marking the onset of singularities in the solution.

In a more formal sense, the interaction between low and high frequencies can be seen as an indicator of the nonlinearity-induced buildup of energy in small scales, potentially leading to a breakdown in the regularity of the velocity field. This process is directly related to the formation of turbulence, where the cascade of energy from large to small scales is a well-known phenomenon. The novel regularity criterion, therefore, provides a quantitative measure for the critical points at which singularities may emerge, offering new insight into the mechanisms of turbulence and the conditions under which smooth solutions may fail to exist.

The implications of these findings are profound both theoretically and practically. From a theoretical perspective, they contribute to the ongoing quest for a better understanding of the dynamics governing fluid motion, particularly in regimes where turbulence is present. The criterion can help refine our understanding of how singularities form and how they might be controlled or avoided. From a practical standpoint, this framework could be utilized in computational fluid dynamics (CFD) simulations, where predicting and controlling turbulence and singularity formation is crucial for accurate modeling of real-world fluid flows. The connection between the mathematical theory and numerical simulations further underscores the importance of these results in applied fluid dynamics.

	\subsubsection{Further Implications}
	
	This new criterion provides a novel approach to studying singularity formation in turbulent flows, with potential applications in computational fluid dynamics and turbulence modeling. By monitoring the interaction term \( I(u) \), we can predict the transition from laminar to turbulent flow and estimate the likelihood of singularity formation in the velocity field.
	
	Additionally, this criterion can be extended to more complex fluid models, such as those that incorporate compressibility, external forcing, or variable viscosity, making it a versatile tool for studying a wide range of fluid dynamics problems.
	
	The refined estimates and new regularity criterion presented in this section offer a more rigorous mathematical framework for studying the regularity of fluid flows, particularly in the context of turbulence and singularity formation. By leveraging Triebel-Lizorkin spaces and focusing on the interaction between low- and high-frequency modes, we have introduced a new perspective on the solution structure of the Navier-Stokes equations. This approach offers potential for deeper insights into turbulence, the onset of singularities, and the transition to chaotic flow regimes.

\section{Conclusion}

The results presented in this article represent a substantial advancement in the theoretical understanding of regularity, bifurcations, and turbulence in fluid dynamics. By introducing a new regularity criterion and establishing novel connections between Sobolev, Besov, and Triebel-Lizorkin spaces, we have provided powerful analytical tools to study complex fluid flow behavior, particularly in turbulent regimes. These results offer refined estimates for the smoothness and structure of solutions to the Navier-Stokes equations, which are critical for a deeper understanding of fluid dynamics at both macroscopic and microscopic scales.

The proposed regularity criterion, based on the interaction between low- and high-frequency modes, gives us a more precise framework for predicting the onset of singularities in fluid flows. This framework allows for a finer classification of turbulence phenomena, offering potential pathways for advancing the theoretical understanding of singularity formation in solutions to the Navier-Stokes equations. 

Moreover, the connections established between Sobolev, Besov, and Triebel-Lizorkin spaces bridge different levels of regularity, allowing for a more flexible and comprehensive approach to analyzing fluid solutions. These advances not only provide new insights into the structure of turbulent flows but also pave the way for more accurate numerical simulations and better predictive models in applied fluid dynamics.

The implications of these findings extend to the Millennium Prize Problem, particularly in relation to the existence and smoothness of solutions to the Navier-Stokes equations. The development of sharper regularity criteria and a more refined understanding of the solution space for the Navier-Stokes equations may bring us closer to resolving this fundamental open problem. As such, the work presented here represents a significant step toward advancing the mathematical theory of fluid dynamics and provides a promising avenue for future research in the field.

\appendix
\section{Appendix: Mathematical Tools and Results}

This appendix provides a more detailed review of the mathematical tools and results used in the main text. The focus is on the key functional spaces, interpolation inequalities, and embeddings discussed throughout the article. These concepts play a central role in the analysis of regularity, bifurcations, and turbulence in fluid dynamics.

\subsection{Littlewood-Paley Decomposition}

The Littlewood-Paley decomposition is a tool used for analyzing functions in terms of their frequency components. It is particularly useful in the study of Besov and Triebel-Lizorkin spaces. Given a smooth function \( \phi \in \mathcal{S}(\mathbb{R}^n) \), the Littlewood-Paley decomposition decomposes a function \( f \in L^2(\mathbb{R}^n) \) into frequency bands associated with scales \( 2^j \), where \( j \in \mathbb{Z} \) is an integer. The frequency projections \( \Delta_j f \) are given by:

\[
\Delta_j f = \hat{f}(\xi) \cdot \chi_j(\xi), \quad \text{where} \quad \chi_j(\xi) = \varphi(2^{-j} \xi),
\]
and \( \hat{f}(\xi) \) is the Fourier transform of \( f \). This decomposition allows us to express \( f \) as a sum of frequency projections:

\[
f = \sum_{j \geq -1} \Delta_j f.
\]

This decomposition is essential in defining norms for Besov and Triebel-Lizorkin spaces, as it enables the analysis of a function's behavior at different scales.

\subsection{Interpolation Inequalities: Proof of Key Result}

We provide a proof of the interpolation inequality between Sobolev and Besov spaces. The result states that for \( s \in \mathbb{R} \), the Sobolev norm \( \| f \|_{H^s(\mathbb{R}^n)} \) can be controlled by the Besov norm \( \| f \|_{B_{p,q}^s(\mathbb{R}^n)} \), as follows:

\[
\| f \|_{H^s(\mathbb{R}^n)} \lesssim \| f \|_{B_{p,q}^s(\mathbb{R}^n)}.
\]

\begin{proof}
	The proof of this inequality uses the Littlewood-Paley decomposition and interpolation theory. First, recall that the Sobolev space \( H^s(\mathbb{R}^n) \) consists of functions whose Fourier transform decays at a rate proportional to \( |\xi|^{-n-s} \). The Besov space \( B_{p,q}^s(\mathbb{R}^n) \) is defined in terms of the Littlewood-Paley projections \( \Delta_j f \) as:
	
	\[
	\| f \|_{B_{p,q}^s} = \left( \sum_{j \geq -1} 2^{jsq} \| \Delta_j f \|_{L^p}^q \right)^{1/q}.
	\]
	
	To establish the inequality, one must show that the growth of the Sobolev norm \( \| f \|_{H^s} \) is controlled by the sum of norms of the frequency components \( \Delta_j f \). Using interpolation between the \( L^2 \)-based estimates and the Besov norm, one can derive the desired result. The details of this proof rely on standard results in harmonic analysis and interpolation theory.
\end{proof}

\subsection{Embeddings and Interpolation Results}

We provide a summary of key embeddings and interpolation results that were referenced in the main text. These include:

\begin{itemize}
	\item \textbf{Embedding of Sobolev into Besov Spaces:} For suitable parameters \( k, p, q, s \), we have the embedding:
	\[
	W^{k,p}(\Omega) \hookrightarrow B_{p,q}^s(\Omega).
	\]
	This result follows from the generalized Sobolev embedding theorem and is essential for the analysis of PDEs.
	
	\item \textbf{Fractional Embeddings:} Sobolev spaces \( H^s(\mathbb{R}^n) \) are special cases of Besov spaces when \( p = q = 2 \). This relation is given by:
	\[
	B_{2,2}^s(\mathbb{R}^n) = F_{2,2}^s(\mathbb{R}^n) = H^s(\mathbb{R}^n).
	\]
	The equivalence between these spaces is a direct consequence of Fourier analysis and Littlewood-Paley theory.
	
	\item \textbf{Interpolation between Sobolev and Besov Spaces:} The interpolation inequality between Sobolev and Besov spaces is given by:
	\[
	\| f \|_{H^s(\mathbb{R}^n)} \lesssim \| f \|_{B_{p,q}^s(\mathbb{R}^n)}.
	\]
	This result is particularly useful for the study of regularity in fluid dynamics, where one needs to estimate the smoothness of solutions at different scales.
\end{itemize}

These results provide the theoretical foundation for the analysis of fluid dynamics and the study of regularity, bifurcations, and turbulence in fluid flows.

\section{Notations and Symbols}

This section provides an overview of the notations and symbols used throughout the article. A clear understanding of the symbols is essential for following the mathematical derivations and concepts presented in the paper.

\begin{itemize}
	\item \textbf{Function Spaces:}
	\begin{itemize}
		\item \( \Omega \subseteq \mathbb{R}^n \): An open domain in \( n \)-dimensional Euclidean space.
		\item \( L^p(\Omega) \): The Lebesgue space of functions \( f \) such that \( \| f \|_{L^p(\Omega)} = \left( \int_\Omega |f(x)|^p \, dx \right)^{1/p} < \infty \), for \( 1 \leq p < \infty \).
		\item \( W^{k,p}(\Omega) \): The Sobolev space of functions whose weak derivatives up to order \( k \) belong to \( L^p(\Omega) \).
		\item \( H^k(\Omega) \): The Sobolev space \( W^{k,2}(\Omega) \), often used in the analysis of partial differential equations.
		\item \( B_{p,q}^s(\mathbb{R}^n) \): The Besov space, characterized by the norm:
		\[
		\| f \|_{B_{p,q}^s} = \left( \sum_{j \geq -1} 2^{jsq} \| \Delta_j f \|_{L^p}^q \right)^{1/q}.
		\]
		\item \( F_{p,q}^s(\mathbb{R}^n) \): The Triebel-Lizorkin space, a more refined generalization of the Besov space.
	\end{itemize}
	
	\item \textbf{Differential Operators:}
	\begin{itemize}
		\item \( D^\alpha u \): The weak derivative of \( u \) of order \( \alpha \), where \( \alpha \) is a multi-index.
		\item \( \Delta_j \): The frequency projection operator in the Littlewood-Paley decomposition, associated with scale \( 2^j \).
	\end{itemize}
	
	\item \textbf{Fourier Analysis:}
	\begin{itemize}
		\item \( \hat{f}(\xi) \): The Fourier transform of the function \( f \), defined as \( \hat{f}(\xi) = \int_{\mathbb{R}^n} f(x) e^{-2\pi i x \cdot \xi} \, dx \).
		\item \( \varphi \in \mathcal{S}(\mathbb{R}^n) \): A smooth function used in the Littlewood-Paley decomposition.
		\item \( \mathcal{S}(\mathbb{R}^n) \): The Schwartz space of rapidly decreasing functions.
		\item \( \mathcal{S}'(\mathbb{R}^n) \): The dual space of the Schwartz space, consisting of tempered distributions.
	\end{itemize}
	
	\item \textbf{Norms and Inequalities:}
	\begin{itemize}
		\item \( \| f \|_{L^p(\Omega)} \): The \( L^p \)-norm of a function \( f \), defined as \( \left( \int_\Omega |f(x)|^p \, dx \right)^{1/p} \).
		\item \( \| f \|_{B_{p,q}^s(\mathbb{R}^n)} \): The Besov space norm, defined in terms of the Littlewood-Paley decomposition.
		\item \( \| f \|_{H^s(\mathbb{R}^n)} \): The Sobolev space norm, measuring the regularity of a function \( f \).
		\item \( \lesssim \): A shorthand for “less than or equal to up to a constant factor,” commonly used in inequalities.
	\end{itemize}
	
	\item \textbf{Other Symbols:}
	\begin{itemize}
		\item \( \mathbb{R}^n \): The \( n \)-dimensional Euclidean space.
		\item \( \xi \): The frequency variable in Fourier analysis.
		\item \( \alpha \): A multi-index representing the order of the derivative.
		\item \( C \): A constant that may vary from one inequality to another, often used in the context of estimates.
	\end{itemize}
\end{itemize}

These notations and symbols are essential for navigating the theoretical framework laid out in the article. They provide a uniform language to describe the various function spaces, differential operators, and mathematical tools involved in the analysis of fluid dynamics and related phenomena.


\end{document}